\RequirePackage{fix-cm}
\documentclass[smallextended]{svjour3}       
\smartqed  
\usepackage{graphicx}
\usepackage{amssymb,amsfonts}
\usepackage{amscd, pdfsync, vpe}
\usepackage{amsmath}

\newtheorem{prop}[theorem]{Proposition}
\newtheorem{coro}[theorem]{Corollary}
\newtheorem{exa}[theorem]{Example}
\newtheorem{defi}[theorem]{Definition}
\newtheorem{rem}[theorem]{Remark}
\newcommand{\bte}{\begin{theorem}\quad  }
\newcommand{\ete}{\end{theorem} }
\newcommand{\bpr}{\begin{prop}\quad  }
\newcommand{\epr}{\end{prop} }
\newcommand{\ble}{\begin{lemma}\quad }
\newcommand{\ele}{\end{lemma}}
\newcommand{\bco}{\begin{coro}\quad }
\newcommand{\eco}{\end{coro} }
\newcommand{\bex}{\begin{exa}\quad \rm }
\newcommand{\eex}{\end{exa} }
\newcommand{\bde}{\begin{defi}\quad \rm }
\newcommand{\ede}{\end{defi} }
\newcommand{\brm}{\begin{rem} \quad \rm}
\newcommand{\erm}{\end{rem} }
\newcommand{\bdm}{\begin{displaymath} }
\newcommand{\edm}{\end{displaymath} }

\newcommand{\lb}{\label}
\newcommand{\lo}{\longrightarrow}

\begin{document}

\title{ On Cogenerating and Finitely cogenerated    $S$-acts
}


\author{Roghaieh Khosravi         \and
        Xingliang Liang \and Mohammad Roueentan 
}


\institute{R. Khosravi \at
              Department of Mathematics, Faculty of Sciences, Fasa University, Fasa, Iran. \\
              \email{khosravi@fasau.ac.ir}           
           \and
          X. Liang\at
              Department of Mathematics, Shaanxi University of Science and Technology,
 Xi$'$an, Shaanxi, P.R. China
 \and
 M. Roueentan \at
 College of Engineering, Lamerd Higher Education Center,
Lamerd, Iran.
}
\date{Received: date / Accepted: date}

\maketitle

\begin{abstract}
In this paper, cogenerating classes of $S$-acts are introduced as those that can be used to cogenerate $S$-acts in an appropriate sense. Next, finitely cogenerated $S$-acts are introduced. Such an $S$-act is then characterized as one whose socle is finitely cogenerated and large in it. Further, we investigate the $S$-acts cogenerating $S_S$, or  generating the injective envelope $E(S)$ of $S_S$.  This leads us to introduce the classes of  cofaithful and subgenerator $S$-acts as the
dual notions of faithful $S$-acts, lying strictly between the classes of generator and faithful $S$-acts. Ultimately, we investigate the relationship between the mentioned classes of $S$-acts.
\end{abstract}

\section{Introduction}

The important, categorical concepts of generating and cogenerating objects play a crucial role in every concrete category. The notion of cogenerator, as the dual concept of generator, is of great importance in category theory. For instance, in most categories, each required injective object can be written as a product of some cogenerator, or more generally, that of a cogenerating set of objects.

Many books have discussed categories and functors also covered the concept of cogenerators. See \cite{kilp2000,HEST73} for example. Moreover, several monographs, including \cite{wis}, have investigated these concepts in the category of modules over rings. In \cite {Nor93}, Normak studied cogenerator $S$-acts in the category of $S$-acts. Also, Knauer and Normak found a relation between cogenerators and subdirectly irreducible $S$-acts in \cite{morita}. Such investigations have been continued later. See \cite{bar, sed}. 
For basic definitions and terminology related to acts over
monoids, we refer the reader to \cite{kilp2000}.

In this paper, we concentrate on the concepts of cogenerating $S$-acts and finitely cogenerated $S$-acts. To do so, we begin with the definition of cogenerators in arbitrary categories. Then, we derive some special properties of cogenerators in the category of $S$-acts.

Let $ \mathcal{C}$ be a concrete category. Recall that an object $C$ of $\mathcal{C}$ is called a \textit{cogenerator} (or a \textit{coseparator}) in $\mathcal{C}$ if the functor $Mor_\mathcal{C}( - , C)$ is faithful. This means that for any $X,Y \in \mathcal{C}$ and any $f,g\in Mor_\mathcal{C}(X, Y)$ with $f\neq g$, there exists $\beta\in Mor_\mathcal{C}(Y, C)$ such that
 $\beta f\neq \beta g$. Equivalently, an object $K$ is a cogenerator in $\mathcal{C}$ if and only if for every $X \in \mathcal{C}$, there exists a monomorphism $X \longrightarrow \prod K$.

Now, we define the categorically generalized notion, namely, cogenerator in arbitrary categories. Let $\mathcal{U}$ be a non-empty set (class) of objects of a category $\mathcal{C}$. An object $B$ in $\mathcal{C}$ is said to be cogenerated by $\mathcal{U}$, or $\mathcal{U}$-cogenerated, if for every pair of distinct morphisms $f,g : A \lo B$ in $\mathcal{C}$, there exists a morphism $h : B \lo U$ with $U \in \mathcal{U}$ and $hf \neq hg$. Then, $\mathcal{U}$ is called a set (class) of cogenerators for $B$.

From now on, we focus on the category of $S$-acts. Let $S$ be a monoid, and $\mathcal{C}$ be a class of $S$-acts.
 An $S$-act $A_S$ is (finitely) cogenerated by $\mathcal{U}$ in case there is a (finite) indexed set $(U_i)_{i\in I}$ in $\mathcal{U}$ and a monomorphism $A\longrightarrow \prod_{i\in I} U_i$.
If $\mathcal{U} =\{U\}$ is a singleton, then we simply say that $U$ (finitely) cogenerates $A$.

In the remainder of this section, we introduce the notion of socle for $S$-acts. To do this, we need to recall some concepts of $S$-acts. An  $S$-act is called \textit{simple} if it contains no subacts other
than  itself, and \textit{ $\theta$-simple}  if it contains no subacts other than
itself and the one element subact $\Theta$. An equivalence relation $\rho$ on an $S$-act $A$ is called
 a \textit{congruence} on $A_S$, if $a \rho a'$ implies $(as)\rho (a's)$
for $a, a' \in A,~ s \in S$.
The set of all congruences on $A_S$ is denoted by $Con(A)$. Clearly, $\Delta_A=\{(a,a)|~a\in A\}, \nabla_A=A\times A \in Con(A)$. Recall from \cite{bar} that a monomorphism $f:A\lo B$ of $S$-acts is said to be \textit{essential} if for each homomorphism $g :B\lo C$, $g$ is a monomorphism whenever $g f$ is.  If $f$ is an inclusion map, then $B$ is said to be an \textit{essential extension} of $A$, or $A$ is called \textit{large} in $B$. In this situation, we write $A\subseteq' B$. It follows from \cite[Lemma 3.1.15]{kilp2000} that $A\subseteq' B$ if and only if for every non-trivial $\theta\in Con(B)$, $\theta \cap \rho_A \neq \Delta_B$.  Recall from \cite{fell} that  if $S$ contains a zero, a non-zero subact $B$ of $A_S$  is called \textit{intersection
large} if for all non-zero subact $ C$  of  $A_S$, $B\cap C \neq \Theta$, and will denoted by $B$  is $\cap$-large in $A_S$.   It follows from  \cite[Proposition 4.7]{fell} that  every large subact of $A_S$  is $\cap$-large, but the converse is not true.

In  module theory, \cite{and},  the socle of a module  is defined to be the sum of the minimal nonzero its submodules.
Equivalently, the intersection of its essential submodules. For $S$-acts, we define two notions as follows.
\bde Let $A$ be a right $S$-acts.
\begin{itemize}
  \item[\rm{(i)}]
$\rm{Soc}(A) = \cap \{L \subseteq A|~L\subseteq' A\}.$
  \item[\rm{(ii)}] If a monoid $S$ contains with a zero, we define
$$S(A) = \cap \{L \subseteq A|~L~ \rm{is} ~\cap\rm{-large}~\rm{ in}~ A\}.$$
\end{itemize} \ede

 If $\rm{Soc}(A)\neq \emptyset$,  then $\rm{Soc}(A)$ is a subact of $A$.
By an argument closely resembles the proof in module theory, one can show that if $S$ contains a zero,
 $$S(A) = \cup \{L \subseteq A|~L~ \rm{is~a ~}~\theta-\rm{simple~subact~of}~ A\}.$$
  Obviously,  $S(A)\subseteq  \rm{Soc}(A)$. But, unlike the case of module theory, the converse can not be valid in general. For instance,  if $S=(\mathbb{N}, max)\cup\{\infty\}$, it is not difficult to see that $S(S)=\{\infty\}\varsubsetneq  \rm{Soc}(S)=S$.

\section{Cogenerating $S$-acts}

In this section, we consider cogenerators for a class of $S$-acts as a common generalization of cogenerators in the category of $S$-acts, and obtain some characterizations of cogenerators.

\bde Let $\mathcal{C}$ be a class of $S$-acts. An $S$-act $A$ is \textit{(finitely) cogenerated} by $\mathcal{C}$ (or
$\mathcal{C}$ (finitely) cogenerates $A$) in case there is a (finite) indexed set $(C_i)_{i\in I}$ in $\mathcal{C}$
and a monomorphism
$A \longrightarrow \prod_{i\in I} C_i$.\ede
We make the obvious adjustments in terminology if $\mathcal{C} =
\{C\}$ is a singleton. 
The class of all $S$-acts cogenerated by $\mathcal{C}$ is denoted by $Cog(\mathcal{C})$. Also,
$FCog(\mathcal{C})$ denotes the class of all $S$-acts which
finitely cogenerated by $\mathcal{C}$.
An $S$-act $C$ is a cogenerator for $Cog(\mathcal{C})$ in case
$Cog(\mathcal{C}) = Cog(C)$. A cogenerator for the class of all right $S$-acts is simply called  a \textit{cogenerator}, without any reference to the class.

The proof of the following proposition is similar to that of \cite[Theorem 2.4.18]{kilp2000}.
\bpr Let $A_S$ and $B_S$ be $S$-acts. The following are equivalent.
\begin{itemize}
	\item[\rm{(i)}]
	$B_S$ (finitely) cogenerates $A_S$.
	\item[\rm{(ii)}] There exists a (finite) subset 	$H$ of $\rm{Hom}(A,B)$ with $\cap_{h\in H}\ker h=\Delta_A$.
	\item[\rm{(iii)}] For any $S$-act $X_S$ and any $f,g\in \rm{Hom}(X, A)$ with $f\neq g$, there exists $\beta\in \rm{Hom}(A, B)$ such that
	$\beta f\neq \beta g$.
\end{itemize} \epr
Now we give a definition of the cotrace of a class of $S$-acts which will be useful to characterize cogenerators later on.
 \bde \lb{de2.30} Let $\mathcal{C}$ be  a class of $S$-acts.
 The \textit{cotrace } of  $\mathcal{C}$ in an $S$-act $A_S$ is defined by\\\centerline{ $cotr_{A_S}(\mathcal{C})= \cap\{\ker g|~g:A_S\longrightarrow C,\ \rm{for~some~} C\in \mathcal{C} \} =\cap_{C\in\mathcal{C}}  Cog(C)$.}\ede

In particular, when $\mathcal{C} = \{C\}$ is a singleton, Definition \ref{de2.30} is the definition of the cotrace of $C$ in $ A$ as mentioned in \cite[Definition 2.4.16]{kilp2000}, which is denoted
by

$$
	cotr_{A_S}(C_S)= \cap_{g\in Hom(A_S,C_S)}  \ker g.
$$
 Note that $cotr_{A_S} (\mathcal{C})$ is a subact of $A_S\prod  A_S$. By \cite[Theorem 2.4.18]{kilp2000}, an $S$-act $C$ is a \textit{cogenerator}
in the category of $S$-acts if and only if for every $S$-act $A_S$, $cotr_{A_S}(C_S) = \Delta_A$.

\bpr Let $\mathcal{C}$ be a class of $S$-acts, and let $A_S$ be an $S$-act. Then, $cotr_{A_S} (\mathcal{C})$  is the unique smallest congruence $\rho$ of $A_S$ such that $A/\rho$ is
cogenerated by $\mathcal{C}$.\epr

\begin{proof} Let $\{C_i\}_{i\in I}$ be an indexed set in $\mathcal{C}$ such that
$$\rho=cotr_{A_S}(\mathcal{C})= \cap\{\ker g_i|~g_i:A_S\longrightarrow C_i,~ i\in I\}.$$
 Define $g=\prod_{i\in I} g_i: A\lo \prod_{i\in I} C_i$ by $g(a)=(g_i(a))_{i\in I}.$ It can be easily checked that $\rho=\ker g.$ Using the homomorphism theorem for $S$-acts, we find a monomorphism $g':A/\rho\lo\prod_{i\in I} C_i $ with $g'([a]_\rho)=g(a)$. Then $A/\rho\in Cog(\mathcal{C})$. Now, suppose that $\sigma$ is a congruence on $A$ such that $A/\sigma\in Cog(\mathcal{C})$. Then, there exists a monomorphism $f:A/\sigma\lo\prod_{j\in J} C_j$.  So, $f\pi :A\lo\prod_{j\in J} C_j$ and $\sigma=\ker f\pi $. This implies $f_j=\pi_j f \pi:A\lo C_j$, and that $\sigma=\ker f\pi=\cap_{j\in J} \ker f_j$. Thus, $\rho=cotr_{A_S}(\mathcal{C})\subseteq \cap_{j\in J} \ker f_j=\sigma$, and the result follows.   \end{proof}

Using the previous proposition, we obtain the following result.
\bco Let $A_S$ be an $S$-act, and $\mathcal{C}$ be a class of $S$-acts. The following hold. \begin{itemize}
	\item[\rm{(i)}]
 $\mathcal{C}$ cogenerates $A_S$ if and only if $cotr_{A_S} (\mathcal{C})=\Delta_A$.
	\item[\rm{(ii)}]  Let $\sigma$ be a congruence on $A_S$.
Then, $\sigma= cotr_{A_S} (\mathcal{C})$ if and only if $\sigma \subseteq cotr_{A_S} (\mathcal{C})$ and $cotr_{A_S/\sigma} (\mathcal{C})=\Delta_{A/\sigma}$. \end{itemize}\eco

\ble Let $\mathcal{C}$ and $\mathcal{D}$  be two classes of $S$-acts.
 If $\mathcal{D}\subseteq Cog(\mathcal{C})$, then
  $Cog(\mathcal{D})\subseteq Cog(\mathcal{C})$ and
   $cotr_{A_S} (\mathcal{C})\subseteq cotr_{A_S} (\mathcal{D})$ for each $S$-act $A_S$.\ele
\begin{proof}  The first part is obvious. To prove the second part, suppose that  $(a,a')\notin cotr_{A_S} (\mathcal{D})$. Then there exists a homomorphism $f:A \lo D$ with $(a,a')\notin \ker f$ for some $D\in \mathcal{D}$, that is, $f(a)\neq f(a')$. Since $D\in Cog(\mathcal{C})$, there exists a
		homomorphism $h:D \lo C $ with $(f(a),f(a'))\notin  \ker h$ for some $C\in \mathcal{C}$. Now, we obtain $hf:A\lo C$
		with $(a,a')\notin \ker hf$. So, $(a,a')\notin cotr_{A_S} (\mathcal{C})$. \end{proof}

 \bpr\lb{pr2.6}  Let $\mathcal{C}$ be a class of $S$-acts.
 \begin{itemize}
	\item[\rm{(i)}]
 If $A\in Cog( \mathcal{C})~ (FCog (\mathcal{C}))$ and $g:A'\longrightarrow A$ is a monomorphism, then $A'\in Cog( \mathcal{C})~(FCog (\mathcal{C}))$.

\item[\rm{(ii)}] If $(A_i)_{i\in I}\in Cog( \mathcal{C}) ~(FCog (\mathcal{C}))$, then $\prod_{i\in I} A_i$
is in $Cog( \mathcal{C})~ (FCog (\mathcal{C}))$. \end{itemize}\epr
\begin{proof} (i). Let $A\in Cog( \mathcal{C})$ and $g:A'\longrightarrow A$ be a monomorphism. Then there exists a monomorphism $f:A \longrightarrow \prod_{i\in I} C_i$, where $C_i\in \mathcal{C}$ for each $i\in I$. So, $fg:A' \longrightarrow \prod_{i\in I} C_i$ is a monomorphism.

(ii).  Let $(A_i)_{i\in I}\in Cog( \mathcal{C})$. Then for each $i\in I$, there exists a monomorphism $f_i:A_i \longrightarrow \prod_{j_i\in J_I} C_{j_i}$, where $C_{j_i}\in \mathcal{C}$ for each $j_i\in J_I$. So,
$$\prod_{i\in I} f_i:\prod_{i\in I} A_i\longrightarrow \prod_{i\in I} (\prod_{j_i\in J_I} C_{j_i})$$  is a monomorphism, and the result follows.
    \end{proof}

   The following results can be proved similar to Proposition \ref{pr2.6}.
\bpr \lb{pr2.8} If $\mathcal{C}$ is the set $\{C_i | i\in I\}$ of $S$-acts, then the following hold.
\begin{itemize}
	\item[\rm{(i)}]  $Cog(\prod_{i\in I} C_i)\subseteq Cog(\mathcal{C})\subseteq Cog(\coprod_{i\in I} C_i)$.
	
		\item[\rm{(ii)}] If $\rm{Hom}(C_i,C_j)\neq \emptyset$ for any $i,j\in I$, then $\prod_{i\in I} C_i$ and $\coprod_{i\in I} C_i$ are cogenerators for $Cog(\mathcal{C})$. \end{itemize}\epr

\bpr  Let $C$
be a cogenerator for $Cog(\mathcal{C})$. Then for each $S$-act $A_S$,  $cotr_{A_S} (C)= cotr_{A_S} (\mathcal{C})$.
 In particular, if $(C_i)_{i \in I}$ is an indexed set of $S$-acts such that $\hom(C_i,C_j)\neq \emptyset$, then
 $cotr_{A_S} (\prod _{i\in I} C_i)= \cap_{i\in I} cotr_{A_S} (C_i)=cotr_{A_S} (\coprod _{i\in I} C_i)$.
 \epr

 \section{Finitely cogenerated $S$-acts}
 In this section, we focus on finitely cogenerated $S$-acts. In \cite{a note}, for a monoid $S$ with zero, an $S$-act $A_S$ is called \textit{finitely cogenerated} provided that for every non-empty collection $\{A_i
~i\in I\}$ of subacts of $A_S$ with $\cap_{i\in I} A_i = \Theta$, there exists a finite subset $J$ of $I$ such that $\cap_{j\in J} A_j = \Theta$. As we know, the importance of congruences is more than subacts in characterizing the structure of $S$-acts. So, on an arbitrary monoid $S$, we define finitely cogenerated $S$-acts based on congruences and cogenerating sets.

 \bde An $S$-act $A_S$ is called \textit{finitely cogenerated} if for every monomorphism $A \overset{f}{\longrightarrow} \prod_{i\in I} A_i$,
 $$ A \overset{f}{\longrightarrow} \prod_{i\in I} A_i \overset{\pi}{\longrightarrow}\prod_{j\in J} A_j$$ is also a monomorphism for some finite subset $J$ of $I$.\ede

  Clearly, if $A_S$ is finitely cogenerated, then every $S$-act that cogenerates $A_S$ finitely cogenerates $A_S$.
\bpr \lb{pr2.1} For any right $S$-act $A_S$, the following are
equivalent.\begin{itemize}
	
	\item[\rm{(i)}] $A_S$ is finitely cogenerated.
	\item[\rm{(ii)}] For every family of homomorphisms $\{f_i: A \lo A_i\}$ in $S$-acts with  $\cap_{i\in I}\ker f_i=\Delta_A$, there is a finite subset  $J$ of $I$ with $\cap_{j\in J}\ker f_j=\Delta_A$.
 \item[\rm{(iii)}] For any family of congruences $\{\rho_i|~ i\in I\}$ on $A_S$, if $\cap_{i\in I}\rho_i=\Delta_A$, then $\cap_{j\in J}\rho_j=\Delta_A$ for some finite subset $J$ of $I$.

\item[\rm{(iv)}] Every subact of $A_S$ is finitely cogenerated.
\end{itemize}
	 \epr
   \begin{proof} (i) $\Rightarrow$ (ii)  Let  $\{f_i: A \lo A_i\}$ be a family of homomorphisms in $S$-acts with  $\cap_{i\in I}\ker f_i=\Delta_A$. Then, $f=\prod_{i\in I} f_i:A\lo \prod_{i\in I} A_i$ is defined by $f(a)=(f_i(a))_{i\in I}$ and $\ker f=\cap_{i\in I} \ker f_i=\Delta_A$. So, $f$ is a monomorphism, and $\pi f: A \longrightarrow\prod_{j\in J} A_j$ is a monomorphism for some finite  subset $J$ of $I$, by our assumption. Thus, $\cap_{j\in J} \ker f_j=\ker \pi f=\Delta_A$.

   The implications  (ii) $\Rightarrow$ (iii) and (iii) $\Rightarrow$ (i) can be proved in a similar way. The implication (iv) $\Rightarrow$ (i) is clear.

   (ii) $\Rightarrow$ (iv) Let $B_S$ be a subact of $A_S$, and $\{\rho_i|~ i\in I\}$ be a family of congruences on $B_S$  such that $\cap_{i\in I}\rho_i=\Delta_B$. It is clear that $\sigma_i=\rho_i\cup \Delta_A$ is also a   congruence on $A_S$, for each $i\in I$. Since $\cap_{i\in I}\sigma_i=\cap_{i\in I}\rho_i\cup\Delta_A=\Delta_A$, by our assumption, $\cap_{j\in J}\sigma_j=\Delta_A$ for some finite subset $J$ of $I$. Thus $\cap_{j\in J}\rho_j\subseteq\Delta_A\cap (B \times B)=\Delta_B$, and the result follows.\end{proof}

   In the following definition, we use Rees congruences instead of subacts to define a weaker notion.
    \bde An $S$-act $A_S$ is called \textit{finitely Rees cogenerated} whenever for any family of Rees congruences $\{\rho_{B_i}|~ i\in I\}$ on $A_S$, if $\cap_{i\in I}\rho_{B_i}=\Delta_A$, then $\cap_{j\in J}\rho_{B_j}=\Delta_A$ for some finite subset $J$ of $I$.\ede

Clearly, $A_S$ is finitely Rees cogenerated if and only if for any family $\{B_i|~ i\in I\}$ of subacts of $A_S$, if $| \cap_{i\in I}B_i| \leq 1$, then $| \cap_{j\in J}B_j| \leq 1$ for some finite subset $J$ of $I$. Also, for a monoid $S$ with zero, this is equivalent to the following statement. If $\cap_{i\in I}B_i=\Theta$, then $\cap_{j\in J}B_j=\Theta$ for some finite subset $J$ of $I$, the same as it is defined in \cite{a note}. Moreover, every subact of a finitely Rees cogenerated $S$-act is finitely Rees cogenerated.






  Note that every ﬁnitely cogenerated  $S$-act is ﬁnitely Rees cogenerated,  but the following example shows that the converse is not true.
 \bex Let $S = (\mathbb{N}, \min)^\varepsilon = (\mathbb{N}, \min) \dot{\cup} \{\varepsilon\}$, where $\varepsilon$ denotes the externally adjoint identity, and denote $\min\{m,n\}$ by $m\ast n$. Then, $K_S = S\setminus \{\varepsilon\}$ is a right ideal of $S$. The subacts of $K_S$ are $1S\subseteq 2S\subseteq 3S\subseteq \ldots$. Hence, $K_S$ is ﬁnitely Rees cogenerated. We claim that $K_S$ is not finitely cogenerated. For each $n\in  K_S$, define $f_n:  K_S\lo K_S$ by $f_n(m)=m\ast n$. It can be easily checked that $\cap_{n\in \mathbb{N}}\ker f_n=\Delta_K$. But, for each finite subset $J$ of $K$, $\cap_{n\in J}\ker f_n\neq \Delta_K$. Therefore, $K$ is not finitely cogenerated. \eex

 Using unique decomposition Theorem ( \cite[Theorem 1.5.10]{kilp2000}) we shall obtain the structure of  finitely (Rees) cogenerated S-acts.

\bpr  Every finitely (Rees) cogenerated $S$-act is a finite coproduct of indecomposable $S$-acts.\epr

\begin{proof}  Suppose that $A_S$ is  finitely cogenerated. As we know, $A_S$ has a unique decomposition into indecomposable subacts $\{A_i|~i\in I\}$, that is, $A=\coprod_{i\in I}A_i$. Let $B_i=\coprod_{j\neq i}A_j$ for each $i\in I$. Then, $B_i$ is a proper subact of $A_i$ and $\cap_{i\in I}\rho_{B_i}=\Delta_A$. Now, since $A_S$  is finitely cogenerated, there exists a finite subset $J$ of $I$ such that $\cap_{j\in J}\rho_{B_j}=\Delta_A$. If $J\neq I$, then $A_i\times A_i\subseteq \cap_{j\in J}\rho_{B_j}$ for $i\in I\setminus J$, which is a contradiction. Therefore, $A=\coprod_{j\in J}A_j$ and we are done.  \end{proof}

 Recall from \cite[Definition 2.5.31]{kilp2000} that an $S$-act is called \textit{completely reducible} if it is a coproduct of simple subacts. Now, the previous proposition allows us to deduce the following corollary.
 \bco  Every finitely (Rees) cogenerated, completely reducible $S$-act is finitely generated.\eco

\bpr \lb{pr3.12} Every non-zero finitely cogenerated $S$-act contains a minimal congruence. In particular, every non-zero finitely Rees cogenerated $S$-act contains a minimal subact. \epr
\begin{proof} Let $A_S$ be a finitely cogenerated $S$-act, and $\mathfrak{A}$ be the set of all non-diagonal congruences on $A_S$. Then, $\nabla_A=(A\times A)\in\mathfrak{A}$ and $\supseteq$ makes $\mathfrak{A}$ into a poset. Let $\{\rho_i|~i\in I\}$ be a chain in $\mathfrak{A}$. If $\cap_{i\in I}\rho_i=\Delta_A$, then by Proposition \ref{pr2.1}, $\cap_{j\in J}\rho_j=\Delta_A$ for some finite subset $J$ of $I$. Since $\{\rho_i|~i\in I\}$ is a chain, $\rho_k=\Delta_A$ for some $k\in I$, which is a contradiction. Thus, by Zorn's lemma, $\mathfrak{A}$ contains a minimal element.

Replacing congruence with Rees congruence in this proof, we obtain the second part. \end{proof}

\bpr  Let $f:A\lo B$ be an essential monomorphism. If $A_S$ is finitely cogenerated, then so is $B$. In Particular, every essential extension (injective envelope) of a finitely cogenerated $S$-act is again finitely cogenerated.
\epr

\begin{proof}   Let $A_S$ be finitely cogenerated, and $f:A\lo B$ be an essential monomorphism. Suppose that $g:B \longrightarrow \prod_{i\in I} A_i$ is a monomorphism. Then, $gf:A \longrightarrow \prod_{i\in I} A_i$ is a monomorphism, and since $A_S$ is finitely cogenerated, $A \overset{gf}{\longrightarrow} \prod_{i\in I} A_i \overset{\pi}{\longrightarrow}\prod_{j\in J} A_j$ is a monomorphism for some finite subset $J$ of $I$. Now, we find that $f$ is an essential monomorphism and $\pi gf$ is a monomorphism. These facts imply that $\pi g$ is a monomorphism. Thus, $B$ is finitely cogenerated.\end{proof}

Now, we consider finitely  (Rees) cogenerated factor $S$-acts.
\bpr\lb{le002} Let $A_S$ be an $S$-act, and $\theta$ be a congruence on $A_S$. Then, $A/\theta$ is finitely (Rees) cogenerated if and only if for any family of (Rees) congruences $\{\rho_i|~ i\in I\}$ on $A_S$, if $\cap_{i\in I}\rho_i=\theta$, then $\cap_{j\in J}\rho_j=\theta$ for some finite subset $J$ of $I$.
\epr

\begin{proof} \textit{Necessity.} Let $\theta$ be a congruence on an $S$-act $A_S$ such that $A/\theta$ is finitely cogenerated. Let $\cap_{i\in I}\rho_i=\theta$, where $\rho_i\in Con(A)$ for each $i\in I$. Define
$$\overline{\rho_i}=\{([a]_\theta,[b]_\theta)|~(a,b)\in \rho_i\}.$$
It can be easily checked that $\overline{\rho_i}\in Con(A/\theta)$ and $\cap_{i\in I}\overline{\rho_i}=\Delta_{A/\theta}$. By our assumption, $\cap_{j\in J}\overline{\rho_j}=\Delta_{A/\theta}$ for some finite subset $J$ of $I$. Thus, $\cap_{j\in j}\rho_j=\theta$.

\textit{Sufficiency.} To show that $A/\theta$ is finitely cogenerated, suppose that $\cap_{i\in I}\sigma_i=\Delta_{A/\theta}$, where $\sigma_i\in Con(A/\theta)$ for each $i\in I$. Define
$$\rho_i=\{(a,b)|~([a]_\theta,[b]_\theta)\in \sigma_i\}.$$
It can be easily checked that $\rho_i\in Con(A)$ and $\cap_{i\in I}\rho_i=\theta$. By our assumption, $\cap_{j\in J}\rho_j=\theta$ for some finite subset $J$ of $I$. Therefore, $\cap_{j\in j}\sigma_j=\Delta_{A/\theta}$, and the result follows.

The case finitely Rees Cogenerated, It suffices to restrict congruences to Rees congruence.  \end{proof}

Recall from \cite{hollow} that $ \mathrm {Rad} (A)$  is the intersection of all maximal subacts of $A_S$.
If $A_S$  contains no maximal subacts, we let $\mathrm {Rad} (A)= A$. If $\mathrm {Rad} (A)\neq \emptyset$, then $\mathrm {Rad} (A)$ is a subact of $A_S$. Now we further consider the factor act  $A/{\mathrm {Rad} (A)}$.
   \bpr    If  $A/{\mathrm {Rad} (A)}$ is finitely Rees cogenerated, then it is cogenerated by finitely many $\theta$-simple $S$-acts. Moreover, $A_S$ has only finitely many maximal subacts.\epr
  	\begin{proof} If $A_S$  contains no maximal subacts, since $\mathrm {Rad} (A)= A$, the result follows. Otherwise,  suppose that  $A/{\mathrm {Rad} (A)}$ is finitely Rees cogenerated. Let
  $${\displaystyle \mathrm {Rad} (A)=\cap _{i\in I}\{M_i\mid M_i{\mbox{ is a maximal subact of A}}\}\,}.$$
  Define $f: A\lo \prod_{i\in I} A/{M_i}$ by $f(a)=([a_i]_{\rho_{M_i}})$. Then, $f$ is an epimorphism such that $\ker f=\rho_{\mathrm {Rad} (A)}=\cap_{i\in I}\rho_{M_i}$.  Using the homomorphism theorem, $\overline{f}:A/{\mathrm {Rad} (A)}\cong \prod_{i\in I} A/{M_i}$. Since $A/{\mathrm {Rad} (A)}$ is finitely Rees cogenerated, we find that $A/\mathrm {Rad} (A)\cong \prod_{j\in J} A/{M_j}$ for a finite subset $J$ of $I$. Moreover, since $M_i$ is maximal, $A/{M_i}$ is $\theta$-simple and the result follows. To show the second part, let $B=\cap_{j\in J} M_j$, and so  $A/\mathrm {Rad} (A)\cong  A/B$. Then  $\mathrm {Rad} (A)=B=\cap_{j\in J} M_j$, and the set of maximal subacts of $A_S$ is finite.
   \end{proof}

 Proposition \ref{pr3.12} together with the fact that $S(A)\subseteq \rm{Soc}(A)$ yields that if $A_S$ is finitely cogenerated,  then Soc$(A)\neq \emptyset$. Now, we use the concepts of essentiality and socle to characterize finitely cogenerated $S$-acts.

\bte\lb{te2.3} An $S$-act $A_S$ is finitely cogenerated if and only if Soc$(A)$ is a finitely cogenerated  large subact of $A_S$.

\ete

 \begin{proof}  If $A_S$ is finitely cogenerated, then the same is also true for each of its subacts, and in particular for Soc($A$). To prove Soc($A$)$\subseteq' A$, suppose that $\theta\in Con(A)$ satisfies $\theta \cap \rho_{\rm{Soc}(A)} \neq \Delta_A$. It is clear that $\rho_{\rm{Soc}(A)}=\cap_{L\subseteq' A} \rho_L$. So, $(\cap_{L\subseteq' A} \rho_L)\cap\theta=\Delta_A$. Since $A_S$ is finitely cogenerated, there exist $L_1,\ldots,L_n\subseteq' A$ such that $(\cap_{i=1}^{i=n}\rho_{L_i})\cap\theta=\Delta_A$. Then, the fact that each $L_i$ is large implies  $\theta=\Delta_A$, and the result follows.

 On the other hand, every essential extension of a finitely cogenerated $S$-act is again finitely cogenerated.  \end{proof}


Recall that  an $S$-act $A_S$ is said to be a \textit{subdirect product} of the family $\{A/{\rho _i}|~i\in I\}$  if $\cap_{i\in I}\rho_i=\Delta_A$. This means that the natural epimorphisms $\pi_i : A_S\lo A/{\rho _i}$ form a monomorphic family. An $S$-act $A_S$ is called \textit{subdirectly irreducible} if every set of congruences $\{\rho_i|~i\in I\}$ on $A_S$ with
$\cap_{i\in I} \rho_i=\Delta_A$ contains $\Delta_A$. Also, an $S$-act $A_S$ is called \textit{irreducible} if
any intersection of a finite number of non-diagonal congruences is non-diagonal.
It is clear that every subdirectly irreducible $S$-act is finitely cogenerated.
The following result can be deduced from the definition of being subdirectly irreducible.
 \bpr \lb{pr3.11} For any right $S$-act $A_S$, the following are
 equivalent.\begin{itemize}
 	
 	\item[\rm{(i)}] $A_S$ is subdirectly irreducible.
 	\item[\rm{(ii)}] There exist distinct elements $a$ and $a'$ of $A_S$ such that every morphism $f:A\lo B$ with $(a,a')\notin \ker f$ is a monomorphism.
 	\item[\rm{(iii)}] There exist distinct elements $a$ and $a'$ of $A_S$ such that $\rho(a,a')$ is the minimum proper congruence of $A_S$.
 	 \item[\rm{(iii)}]	If $f:A \longrightarrow \prod_{i\in I} A_i$ is a monomorphism, then $\pi f:A \longrightarrow  A_j$  is already a monomorphism for some $j\in I$.
 	\item[\rm{(iv)}]  Every subact of $A_S$  is subdirectly irreducible.
 	\item[\rm{(v)}] $A_S$  is a finitely cogenerated irreducible $S$-act.
 \end{itemize}
 \epr

   By Birkhoff's theorem for acts, \cite{gra}, any non-trivial $S$-act is a subdirect product of subdirectly irreducible $S$-acts. Now, using part (iii) of Proposition \ref{pr2.1}, we obtain the following result.
 \bco If $A_S$ is a finitely cogenerated $S$-act, then it is isomorphic to a subdirect product of finitely many subdirectly irreducible $S$-acts. \eco



  \section{Characterization of  the $S$-acts cogenerating  $S$ }

From \cite[Proposition 2.6]{tor}  it follows that  $ S_S$ cogenerates an $S$-act  $A_S$   if and only if  $A_S$  is torsionless, such acts are characterized in \cite{tor}.  Let us turn to the question when an act cogenerates $S$.
 This section  concerns with the properties of $S$-acts which (finitely) cogenerate $S_S$,  or generating the injective envelope $E(S)$ of $S_S$. We introduce the classes of (strongly) cofaithful $S$-acts and we give characterizations of monoids $S$ such that all faithful acts are cofaithful.

  Let $A_S$ be an $S$-act and $a \in A_S$. Then, $\lambda_a:S_S\lo A_S$ is defined by $\lambda_a(s) = as$ for every $s\in S$. The kernel congruence $\ker \lambda_a$ on $S_S$ is called the \textit{annihilator congruence} of $a\in A_S$.
   Recall form \cite{fully stable} that the \textit{right annihilator} of $A_S$ is defined by
$$R_S(A) =\{(s, t)\in S\times S|~as=at,\ \rm{for ~all~} a\in A\},$$
    which is a two-sided congruence on $S$.
   We call $A_S$ a \textit{faithful right} $S$-act if for $s,t\in S$, the equality $as= at$ for all $a\in A$ implies $s= t$. Clearly, $R_S(A) =\cap_{a\in A} \ker \lambda_a$ and $A_S$ is faithful in case $R_S(A)=\Delta_S$.
   On the other hand, for each right $S$-act $A_S$,
$$cotr_S(A_S)=\cap_{g\in Hom(S_S,A_S)}  \ker g=\cap_{a\in A}  \ker \lambda_a=R_S(A).$$

The next theorem will be a  useful description of faithful $S$-acts.
 \bte \lb{te2.10}  For each right $S$-act $A_S$,  the following are equivalent.\begin{itemize}
\item[\rm{(i)}]  $A_S$ is faithful.
	\item[\rm{(ii)}]  $A_S$ cogenerates $S$.
 \item[\rm{(iii)}]
 $A_S$ cogenerates every projective  $S$-act.
  \item[\rm{(iv)}]
 $A_S$  cogenerates every free $S$-act.
 \item[\rm{(v)}]
 $A_S$ cogenerates a generator $S$-act. \end{itemize}
   \ete
 \begin{proof} Since $cotr_S(A_S)=R_S(A)$, clearly (i) and (ii) are equivalent. It suffices to show   (ii) $\Rightarrow $  (iii).
 Let $S\hookrightarrow A^J$. Suppose that $P=\coprod_{i\in I} e_iS$  is a projective $S$-act. Since $e_iS$  is a retract of $S$ and $S\in Cog(A)$, we deduce that $e_iS\in Cog(A)$. Using Proposition \ref{pr2.8}, since $\hom(A,A)\neq \emptyset$, $\coprod_{i\in I} A$ is a cogenerator for $Cog(A)$.  If $f_i: e_iS\rightarrow A^J$, then $\coprod_{i\in I} f_i:\coprod_{i\in I} e_iS\hookrightarrow \coprod_{i\in I} A^J\hookrightarrow (\coprod_{i\in I} A)^J$. Therefore, $P\in Cog(\coprod_{i\in I} A)=Cog(A)$, as desired. \end{proof}

From the previous theorem, we know that faithful $S$-acts can be characterized as those $S$-acts cogenerate
 $S_S$ or, equivalently, cogenerate every projective $S$-act.
In the category of  modules, the concept of  cofaithful is the dual notion of faithful as the modules which generate every injective module, which is equivalent to modules finitely cogenerate  $R$, such modules
are also called  subgenerators of Mod-$R$.  Unlike the case  for modules, these properties  are no longer valid  for $S$-acts.
  This description allows us to define the following dual notions:
\bde Let $A_S$ be an $S$-act. \begin{itemize}
\item[\rm{(i)}]  $A_S$  is called \textit{cofaithful} in case  $A_S$  finitely cogenerates  $S_S$, i.e., there exists a positive integer $n$  such that $S_S$ can be embedded to $A^n$.
\item[\rm{(ii)}] $A_S$  is called \textit{subgenerator} in case it generates every injective $S$-act. \end{itemize}\ede

\ble \lb{le4.1} Let $A_S$  be a right $S$-act. The following are equivalent: \begin{itemize}
	\item[\rm{(i)}]  $A_S$ is cofaithful.
\item[\rm{(ii)}] There exists a finite subset $B$ of elements of $A_S$  such that
$R_S(B) = \Delta_S$.
\end{itemize}\ele

\begin{proof}  (i) $\Rightarrow$ (ii).  Let $f: S\rightarrow A^n$ be a monomorphism. Suppose $f(1)=(a_1,...,a_n)$, and set $B=\{a_1,...,a_n\}$. If $(s,t)\in R_S(B)$, then $a_is=a_it$ for each $1\leq i\leq n$. Clearly, $f(s)=f(t)$, and thus $s=t$.

 (ii) $\Rightarrow$ (i). If $R_S(\{a_1,..., a_n\}) = \Delta_S$, then $\lambda_{(a_1,...,a_n) } :S\rightarrow A^n$  is a monomorphism.
 \end{proof}

\bco \lb{co4.4}  A cofaithful $S$-act contains a finitely generated faithful
subact. Moreover, if $S$ is a commutative monoid, the converse is valid.\eco

\begin{proof}  Let $A_S$  be  cofaithful. Then there exists a finite subset $\{a_1,...,a_n\}$  of  $A_S$  with $R_S(\{a_1,..., a_n\}) = \Delta_S$. Set $B=\cup_{i=1}^{i=n} a_iS$. Clearly,  $R_S(B) \subseteq R_S(\{a_1,..., a_n\}) =\Delta_S$, and $B$ is faithful.
To show the second part, suppose that a finitely generated subact $B=\cup_{i=1}^{i=n} a_iS$ of $A_S$  is faithful. Since $S$ is commutative, $ R_S(\{a_1,..., a_n\}) =R_S(B)=\Delta_S$, and thus $A_S$ is cofaithful.
 \end{proof}
 As we know, the notion $A_S$  finitely cogenerates  $S_S$  means there exists a positive integer $n$  such that $S_S$ can be embedded to $A^n$. If $n=1$, i.e., $S\hookrightarrow A$, we say that $A_S$  \textit{cyclically cogenerates}  $S_S$.
 \bpr  Let $A_S$  be a right $S$-act. The following are equivalent: \begin{itemize}
	\item[\rm{(i)}]  $A_S$ is a  subgenerator.
\item[\rm{(ii)}] $A_S$  generates  $E(S)$.
\item[\rm{(iii)}] $A_S$  cyclically cogenerates  $S_S$.
\item[\rm{(iv)}] There exists an element $a\in A$  such that $R_S(\{a\}) = \ker \lambda_a=\Delta_S$.
\item[\rm{(v)}] $A_S$  contains a cyclic generator subact.
\end{itemize}\epr

\begin{proof}   By an argument similar to that of Lemma \ref{le4.1}, one can prove  (iii) $\Leftrightarrow$  (iv). The implications  (iv) $\Leftrightarrow$  (v) and (i) $\Rightarrow$  (ii) are clear.

(ii) $\Rightarrow$  (iv).  Let $f:\coprod_{I} A\rightarrow E(S)$  be an epimorphism and $\iota:S\hookrightarrow E(S)$. For $\iota(1)\in E(S) $,  $f(a)=\iota(1)$ for some $a\in A$.  It is easy to see that $R_S(\{a\}) =\Delta_S$.

(iii) $\Rightarrow$  (i). Let $E$  be  an injective $S$-act and $f:S\hookrightarrow A$ be an monomorphism. For each $b\in E$, there exists a homomorphism  $g_b:A\rightarrow E$ such that $g_bf=\lambda_b$. So we have the homomorphism $g=\coprod_{b\in E}g_b: \coprod_{b\in E}A\rightarrow \coprod_{b\in E}E$. Hence $\coprod_{b\in E}A\rightarrow \coprod_{b\in E}E\rightarrow E$  is an epimorphism, as desired.
\end{proof}

It is easily checked that  the following  implications are valid,\\\centerline{generator  $\Longrightarrow$  subgenerator $\Longrightarrow$ cofaithful $\Longrightarrow$  faithful.}
The following example shows that these implications are strict.
\bex
 \begin{itemize}
	\item[\rm{(i)}] The implication cofaithful $\Longrightarrow$  faithful is strict:\\
 Let $S =  (\mathbb{N}, \min) \dot{\cup} \{\varepsilon\}$, where $\varepsilon$ denotes the externally adjoint identity, and $A_S = S\setminus \{\varepsilon\}$. Obviously, $A_S$ is faithful but not cofaithful.
\item[\rm{(ii)}] The implication  subgenerator $\Longrightarrow$  cofaithful  is strict:\\
Let  $S =\{1,0, e, f\}$  be the semilattice where $ef = fe = 0$,
and take $A_S=\{e,f,0\}$. Clearly, $A_S$  is not a subgenerator. But $R_S(\{e,f\})=\Delta_S$, and so $A_S$ is cofaithful.

\item[\rm{(iii)}]  The implication  generator $\Longrightarrow$  subgenerator  is strict:\\
Let $S =  (\mathbb{N}, .) $, and $A_S=\mathbb{N}\coprod^{2\mathbb{N}} \mathbb{N}$. By \cite[Example 2.2.]{sed}, $A_S$ is not a generator. But $S\hookrightarrow A$, and so $A_S$ is a subgenerator.

\end{itemize}
\eex
Concluding this section, we pointed at the conditions on a monoid $S$ under which the converses of implications are true. Recall that a monoid $S$ is said to be \textit{right  self-injective}  if  $S_S$  is injective.
\bte For a monoid $S$  the following are equivalent: \begin{itemize}
	\item[\rm{(i)}] Every subgenerator $S$-act is a generator.
	\item[\rm{(ii)}] $E(S)$  is a generator.
	\item[\rm{(iii)}]  $S_S$ is right self-injective.\end{itemize}\ete

 \begin{proof}   (i) $\Rightarrow$  (ii) is clear.

   (ii) $\Rightarrow$  (iii). Since $E(S)$  is a generator, $S$ is a retract of $E(S)$,  and so $S_S$ is injective.

  (iii) $\Rightarrow$  (i). Let $A_S$  be  a subgenerator $S$-act. Then $f:S\hookrightarrow A$, and injectivity of $S_S$ implies that $S_S$ is a retract of $A_S$. Thus $A_S$ is a generator.

 \end{proof}
By an argument similar to that of Corollary \ref{co4.4}, It is not difficult to obtain the following result.
\bpr If $S_S $  is  irreducible, then every cofaithful $S$-act is a subgenerator.  Moreover, if $S$ is a commutative monoid, the converse is true.  \epr

 Using  Theorem \ref{te2.10}, the following  result can be obtained.

\bpr $S_S $  is  finitely cogenerated if  and only if every faithful $S$-act is cofaithful.  \epr

\end{document}